\documentclass[12pt]{amsart}

\usepackage{url}
\usepackage{epsfig}
\usepackage{amsthm}
\usepackage{amsmath}
\usepackage{amssymb}
\usepackage{amsfonts} 

\theoremstyle{plain}
\newtheorem{theorem}{Theorem}[section]
\newtheorem{lemma}[theorem]{Lemma}
\newtheorem{corollary}[theorem]{Corollary}
\newtheorem{problem}[theorem]{Problem}
\newtheorem{remark}[theorem]{Remark}
\newtheorem{proposition}[theorem]{Proposition}
\newtheorem{example}[theorem]{Example}

\theoremstyle{definition}
\newtheorem{definition}[theorem]{Definition}

\newcommand{\Z}{\mathbb{Z}}
\newcommand{\R}{\mathbb{R}}

\newcommand{\N}{\mathbb{N}}

\newcommand{\cone}{\mathrm{cone}}

\begin{document}
\title{{\bf Supernormal Vector Configurations}}
\author{ Serkan Ho\c{s}ten}
\address{Department of Mathematics \\
     San Francisco State University \\
     San Francisco, CA 94132}
\email{\tt serkan@math.sfsu.edu}
\author{Diane Maclagan}
\address{School of Mathematics \\
        Institute for Advanced Study \\
        Princeton, NJ 08540}
\email{\tt maclagan@ias.edu}
\author{Bernd Sturmfels}
\address{Department of Mathematics\\
 University of California\\
Berkeley, CA 94720.}
\email{bernd@math.berkeley.edu}

\thanks{The second author was supported by NSF grant DMS 97-29992
through IAS.
The third author was supported by NSF grant
DMS 99-70254.}


\begin{abstract}
A configuration of lattice vectors is supernormal
if it contains a Hilbert basis for every cone spanned
by a subset. We study such configurations from 
various perspectives, including triangulations, integer programming
and Gr\"obner bases. Our main result is a bijection between
virtual chambers of the configuration and virtual initial ideals of the
associated binomial ideal.

\end{abstract}

\maketitle

\section{Introduction}

Let $ B = \{b_1,...,b_n \} \subseteq \Z^m$ and let $\cone(B)$ be the
polyhedral cone in $\R^m$ spanned by $B$.  The configuration $B$ is
{\em normal} if every lattice point in $\cone(B)$ is a non-negative
integer combination of $B$.
We say that $B$ is {\em supernormal } if, for every subset $B'$ of
$B$, every lattice point in $\, \cone(B') \,$ is a non-negative
integer combination of $\, B \, \cap \, \cone (B')$.

In Section 2 we discuss supernormal configurations
in low dimensions.  In particular, we
exhibit a finitely generated submonoid of $\mathbb Z^3$ which cannot
be generated by a finite supernormal subset.  This implies that in
general the process of normalization \cite[Algorithm
13.2]{GB+CP} cannot be extended to produce a finite supernormal
generating set.

In Section 3 we characterize  supernormal vector configurations
in terms of polyhedral geometry (triangulations)
and  in terms of integer programming (total dual integrality).
This will generalize the familiar characterizations of
unimodular configurations \cite[\S 8]{GB+CP}. Recall
that a configuration $B$ in $\Z^m$ is {\em unimodular} if, for every
subset $B'$ of $B$, every lattice point in $\cone(B')$ is a
non-negative integer combination of $B'$.

The algebraic theory of integer programming is closely related to
Gr\"obner bases of binomial ideals \cite{ST}.
We encode our configuration $B$ as the  ideal $J_B$
in the polynomial ring $k[x_1,\ldots,x_n]$ generated by
\begin{equation}
\label{Ubino}
\prod_{i: u \cdot b_i > 0} x_i^{u \cdot b_i} \,\,\, - \,
\prod_{j: u \cdot b_j < 0} x_j^{- u \cdot b_j}
\qquad
\hbox{where $u$ runs over $\Z^m$.}
\end{equation}
In the language of  \cite{HT} or \cite[\S 3.3]{SST},
the ideal $J_B$ is the {\it lattice ideal} for the lattice spanned
by the rows of the $(m \times n)$-matrix 
$(b_1,\ldots,b_n)$.

Every vector $w \in \cone(B)$ defines an {\em initial ideal} of
$J_B$ as follows: $in_w(J_B)$ is generated by the monomials $\,
\prod_{i: u \cdot b_i > 0} x_i^{u \cdot b_i} \,$ where $u \in \Z^m$
satisfies $u \cdot w > 0$ and the binomials of the form in equation
(\ref{Ubino}) where $u \in \Z^m$ satisfies $u \cdot w = 0$.  Two
vectors $\, w , w' \in \cone(B)$ lie in the same cell of the {\em
Gr\"obner fan of $J_B$} if $in_w(J_B) = in_{w'}(J_B)$, and they lie in
the same cell of the {\em chamber complex of $B$} if, for every subset
$B' $ of $B$, $w \in \cone(B')$ if and only if $w' \in \cone(B')$.  In
Section 4 we prove:

\begin{theorem} \label{gb=chamber}
If the configuration $B$ is supernormal  then the
chamber complex of $B$ coincides with the
Gr\"obner fan of $J_B$.
\end{theorem}

We note that the converse statement does not hold, even for $m=1$.  
For the special case when
$B$ is unimodular, this theorem follows from \cite[Proposition
8.15]{GB+CP} via Gale duality.  Our proof will be self-contained.

\begin{figure}
\epsfig{file=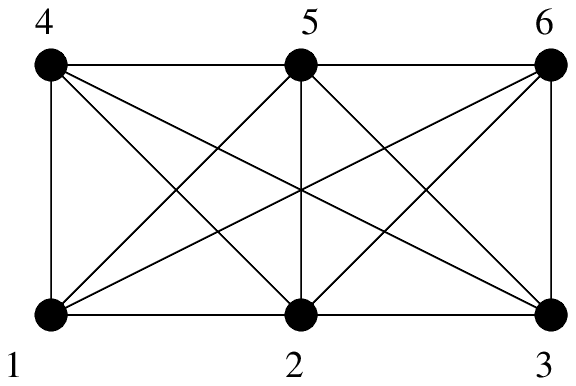, height={4cm}}
\caption{Chamber complex of a rectangle}
\label{2x3grid}
\end{figure}

A longstanding conjecture \cite{ST} states that
the number of facets of any chamber in the Gr\"obner fan of
$J_B$ is bounded by a function of $m$
alone, independent of the coordinates of the $b_i$.
In Section 5 we examine this question for  the supernormal configuration
\begin{equation}
\label{Bintheplane}
 B \quad = \quad \bigl\{\,
(1,u,v) \in \Z^3 \,\, : \,\, (u,v) \in P \cap \Z^2 \,\bigr\} 
\end{equation}
associated with a convex lattice polygon $P$ in the plane.  The
chamber complex of $B$ is
gotten by drawing the line segments connecting any two lattice
points in $P$ as in Figure \ref{2x3grid}.  
It is an open question whether polygons with
arbitrarily many edges can appear in such a picture. 
See Proposition \ref{world-record} for the current status of the 
problem.

The chambers of a vector configuration $B$ are in bijection with the
{\em regular triangulations} of a Gale dual configuration $A$. This
was extended in \cite{DHSS} to a bijection between all
triangulations of $A$ and {\em virtual chambers} of $B$. We reexamine
these concepts in Section 6, and we introduce the following algebraic
analogue: A monomial ideal $M$ in $k[x_1,\ldots,x_n]$ is a {\em
virtual initial ideal} of $J_B$ if $M$ has the same Hilbert function
as $J_B$ with respect to the finest grading which makes $J_B$
homogeneous.  In \cite[\S 10]{GB+CP} such $M$ were called {\em
$A$-graded monomial ideals}. There is a map from virtual initial
ideals to virtual chambers (defined by \cite[Theorem
10.10]{GB+CP}) but this map is in general neither injective nor
surjective.  Our main result is the following extension of Theorem
\ref{gb=chamber}.

\begin{theorem} \label{unimodular2}
If the configuration $B$ is supernormal then the
map from virtual initial ideals of $J_B$ to virtual chambers
of $B$ is a bijection.
\end{theorem}

\section{Examples and Counterexamples} \label{egs}

In this section we study examples of supernormal configurations in low 
dimensions. Recall that a configuration $B$ of vectors in $\Z^m$ is {\em
normal} if it generates the monoid $\Z^m \, \cap \, \cone(B)$.
We call $B$ {\it pointed} if there
exists $u \in \R^m$ such that $b_i \cdot u > 0$ for all $i$. We say
that $B$ is a {\it Hilbert basis} if $B$ is pointed and minimally
generates the monoid $\Z^m \, \cap \cone(B)$.  Clearly, if $B$ is a
Hilbert basis then $B$ is normal.

\vskip .2cm

\noindent {\bf Dimension one: } If $m = 1$ then $B$ is a set of
integers $\, \{ b_1 < b_2 < \cdots < b_n \} \,$ which we assume to be
ordered, distinct, and not containing $0$.  The configuration $B$ is
normal if and only if either $\,b_1 = 1\,$, or $\,b_n = -1\,$, or $b_1
< 0 $, $b_n > 0$ and $\gcd(b_1,\ldots,b_n)= 1$.  But $B$ is
supernormal if and only if either $\,b_1 = 1\,$, or $\,b_n = -1\,$, or
$\{-1,+1\} \subseteq B$.  Thus $B = \{-2,3\} \,$ is normal but not
supernormal. 

The chamber complex of $B$ consists of either one
or two cones, and it coincides with the
Gr\"obner fan of the principal ideal  $J_B$.
For instance, for $B = \{-2,3\} $ the ideal
$\,J_B = \langle x^2 - y^3 \rangle\,$ has two
initial ideals, but for $B = \{2,3\}$ we get
$\,J_B = \langle x^2 y^3 - 1 \rangle\,$ which has
only one initial ideal. This shows that
the converse to Theorem \ref{gb=chamber} does not hold.

The configuration $B$ is pointed if and only if either
$b_1 > 0$ or $b_n < 0$.  Note that in this case $B$ is
normal if and only if either $\{+1\} \subseteq B$
or $\{-1\} \subseteq B$. We conclude that
a one dimensional pointed configuration is
normal if and only if it is supernormal.

\vskip .2cm

\noindent {\bf Dimension two: }
The configuration $B$ consists of distinct nonzero vectors
in the plane $\Z^2$. We assume their ordering
$b_1,b_2,\ldots,b_n$ is counterclockwise,
and we set $ b_{n+1} = b_1$.
They lie in an open half-plane if and only if $B$ is pointed.
The last statement from the $m=1$ case does not hold for $m=2$:
the configuration $B = \{ (1,0), (1,2), (0,1) \} \,$
is pointed and normal but not supernormal.

\begin{proposition} 
A configuration
$B \, = \, \{b_1,\ldots,b_n\} \subseteq \Z^2 \,$ is supernormal if and
only if $\,\, det(b_{i}, b_{i+1}) = 1 \,$ for $i = 1,2,\ldots,n-1 $,
and $det(b_n,b_1)=1$ if $B$ positively spans $\R^2$.
\end{proposition}

\begin{proof}
Suppose $B$ is supernormal. Note that $\cone(b_i,b_{i+1})$ 
contains no other $b_j$, so $b_i$ and $b_{i+1}$ must be a Hilbert basis 
for $\cone(b_i,b_{i+1}) \cap \Z^2$, and thus we have $\det(b_i,b_{i+1})=1$.

Conversely, suppose that $\det(b_i,b_{i+1})=1$ for all $i$.
This means that any lattice point in $\cone(b_i,b_{i+1})$ can be
written as a non-negative integer combination of $b_i$ and $b_{i+1}$.  
Every cone generated by a subset $B'$ of the $b_i$ can be decomposed 
as a union of cones of the form $\cone(b_i,b_{i+1})$, so our assumption 
implies that every lattice point in $\cone(B')$ can be written as a 
non-negative integer combination of the vectors in 
$\,B \cap \cone(B')$. Therefore $B$ is 
supernormal.
\end{proof}

\begin{corollary}
\label{TwoHilb}
Every two-dimensional Hilbert basis is supernormal.
\qed
\end{corollary}

In the language of algebraic geometry, this says that  $B \subseteq \Z^2$
is supernormal if and only if the toric surface $X_\Delta$ is smooth, where
$\Delta$ is the fan whose rays are  the vectors in $B$.  
In higher dimensions, supernormality means that all toric varieties 
that share a fixed Cox homogeneous coordinate ring are smooth.  
This follows from Proposition \ref{uni-triang} below.

\vskip .2cm

\noindent {\bf Dimension three: }
Corollary \ref{TwoHilb} does not hold for $m=3$. Take
$$\,
B \, = \,\bigl\{
(1,0,0),
(0,1,0),
(1,1,1),
(1,1,2),
(1,2,3),
(1,2,4)
\bigr\} \quad \subseteq \quad \Z^3. $$
This is the Hilbert basis for the cone spanned
by $(1,0,0),  (0,1,0)$ and $\, (1,2,4)$.
The configuration $B$ is not supernormal.
To see this consider $B' = \{(0,1,0)$, $ (1,1,1), (1,2,3)\}$
and note that $\,(1,2,2) \,$ lies in
$\cone(B') \,\cap\, \Z^3 \,$ but
not in the monoid generated by
$\,\cone(B') \,\cap\, B\,=\, B'$.
If we add the vector $(1,2,2)$ to $B$ then the resulting
configuration of seven vectors is supernormal.

It is well-known that the monoid of lattice points in any
rational polyhedral cone  has a finite Hilbert basis.
In the previous example,
the Hilbert basis can be enlarged to a finite
supernormal generating set. This raises the
question of whether every rational submonoid of $\Z^m$ is generated
by a finite supernormal subset. This is not the case.

\begin{theorem}
The monoid of lattice points in the three-dimensional cone spanned
by  $P_0 = (-1,1,2),  P_1 = (1,-1,1), P_2 = (0,1,0) $
and $P_3 = (1,0,0)$ is not generated by a finite supernormal subset.
\end{theorem}

\begin{proof}
Consider the following sequence of vectors in this monoid:
$$
P_i  \quad :=  \quad  \frac{1}{2} \cdot
(\, P_{i-2} \,+\, P_{i-1}
\,+ \, P_{i'}) \qquad \hbox{for} \quad i \geq 4,
$$
where $i'= (i \mod 2)$.   Explicitly,
$$ P_{2i}=(0,1,i-1), \quad P_{2i+1}=(1,0,i-1)
\qquad \hbox{for} \quad i \geq 1 $$
At each stage in this iteration,
the three vectors $\, P_{i-2} , P_{i-1} , P_{i'} \,$
generate an index two sublattice of $\Z^3$, and
$P_i$ is the unique vector which completes the Hilbert basis
for their triangular cone.
Suppose there is a finite supernormal generating set $B$ for the ambient
monoid and consider the smallest index $i$ such that $P_i$ is not in $B$.
Then the subset $B' = \{P_{i-2}, P_{i-1}, P_{i'}\}\,$ 
violates  the defining property of $B$ being supernormal.
\end{proof}

While this result shows that not every configuration can be embedded
into a supernormal one, there do exist interesting specific
supernormal configurations in higher dimensions, beyond the familiar
class of unimodular configurations. Here is an example for $m=3$:

\begin{example}
The configuration $B = \{-1,0,+1\}^3$ of
all $27$ vectors whose coordinates have absolute value at most one
is  supernormal.
\end{example}

A  configuration obtained from all lattice points in a 
lattice polytope as in
(\ref{Bintheplane}) is called {\it convex}.  The three dimensional
convex configurations (\ref{Bintheplane}) arising from polygons
play a special role and are discussed in detail in Section
\ref{polygon}.  In Proposition \ref{polygonsupernormal} we show that
they are supernormal.

\vskip .2cm

\noindent {\bf Dimension four and beyond: } 
Most convex configurations in higher dimensions are not supernormal,
however.  Consider the cone over the three-dimensional cube given by the
columns of
$$
 \left( \begin{array}{cccccccc}
1 & 1 & 1 & 1 & 1 & 1 & 1 & 1 \\
0 & 1 & 0 & 1 & 0 & 1 & 0 & 1 \\
0 & 0 & 1 & 1 & 0 & 0 & 1 & 1 \\
0 & 0 & 0 & 0 & 1 & 1 & 1 & 1
\end{array}  \right).
$$
This configuration of eight vectors in $\Z^4$ is convex but not supernormal.
What is missing is the vector  $\,(2,1,1,1) \,$ which represents
the centroid of the cube. The configuration together with
$(2,1,1,1)$ is supernormal.

It would be interesting to identify infinite families
of configurations in higher dimensions
which are supernormal but not unimodular.
Such families might arise from graph theory or
combinatorial optimization.

\section{Polyhedral Characterizations}

In this section we present two characterizations of
supernormal configurations $B$. The first is in terms
of triangulations, and the second involves the
concept of total dual integrality from  integer programming.

A {\it subdivision} of a vector 
configuration  $B = \{b_1,\ldots,b_n\} \subseteq \Z^m $ 
is a polyhedral fan $\Delta$ in $\R^m$ whose support is $\,
\cone(B) \,$ and each of whose rays is spanned by a vector $b_i$
\cite[\S 9]{Zi}.  It is customary to identify $\Delta$
with the collection of subsets $\sigma$ of $B$ which lie
in the maximal cones of $\Delta$.
A subdivision $\Delta$ is {\it regular} if there exists
a vector $c \in \Z^n$ such that $\sigma \subseteq B$ is a face of $\Delta$
if and only if $b_i \cdot x = c_i \text{ for all } b_i \in \sigma$
and $b_i \cdot x < c_i$ otherwise. A subdivision $\Delta$
of $B$ is a {\it triangulation} if each maximal cell 
$\sigma$ has precisely $m$ elements.  A triangulation $\Delta$
of $B$ is {\it unimodular} if every maximal cell  $\sigma$
is a lattice basis for $\Z^m$. 
The triangulation $\Delta$ 
{\it uses all vectors} if each element $b_i$ of $B$
spans a ray of the fan $\Delta$. 

A configuration $B$ is unimodular if and only if every
triangulation of $B$ is unimodular. 
Here it suffices to consider regular triangulations.
We prove an analogous characterization for  supernormal
configurations.

\begin{proposition} \label{uni-triang}
For a configuration $B$,
the following are equivalent:
\begin{enumerate}
\item $B$ is supernormal.
\item Every triangulation of $B$ that uses all vectors is unimodular.
\item Every regular triangulation of $B$ that uses all vectors is
unimodular.
\end{enumerate}
\end{proposition}

\vskip -.3cm

\begin{proof}
We first prove $(1) \Rightarrow (2)$.  Let $B$ be supernormal,
$\Delta$ a triangulation that uses all vectors, and $\sigma =
\{b_{i_1},\ldots,b_{i_m}\}$ a maximal cell of $\Delta$.  If $\sigma$
is not a lattice basis of $\Z^m$ then $\sigma$ does not generate the
monoid $\, \Z^m \, \cap \, \cone(\sigma)$.  Supernormality implies
that $\cone(\sigma)$ contains at least one other vector $b_j \in B
\backslash \sigma$, but then this vector $b_j$ cannot be used in the
triangulation $\Delta$. This contradicts our hypothesis.

The implication $(2) \Rightarrow (3)$ is trivial. 
It remains to show  $(3) \Rightarrow (1)$.
Suppose (3) holds.   Let $B'$ be any subset of $B$.
We construct a regular subdivision of $B$ which has
$\,\sigma = B \, \cap \, \cone(B') \,$
 as one of its faces, and which uses all vectors in
$B \backslash \sigma$ as rays. This subdivision can be 
refined to a regular triangulation $\Delta$ of $B$ that uses all vectors.
By hypothesis, $\Delta$ is unimodular,
and its restriction to $\sigma$ is a unimodular triangulation of $\sigma$.
This implies that $\sigma$ generates the monoid
$\, \cone(B') \, \cap \Z^m$.
We conclude that  $B$ is supernormal.
\end{proof}

Regular subdivisions are polar to the polyhedra with facet normals
in $B = \{b_1,\ldots,b_n\}$. More precisely, for $c \in \Z^n$ we define
the polyhedron
$$ P_c \,\, = \,\, \bigl\{ \, x \in \R^m \,\, : \,\, b_i \cdot x \leq c_i
\,\text{ for } \, i = 1,\ldots,n \,\bigr\}. $$
Let $\mathcal{N}(P_c)$ denote the normal fan of the
convex polyhedron $P_c$.

\begin{lemma} \label{regular=fan}
The normal fan $\mathcal{N}(P_c)$ is a regular
subdivision of $B$. Every regular subdivision of $B$
is the normal fan of $P_c$ for some $c \in \Z^n$.
\end{lemma}

\begin{proof} These statements follow from the fact that $\mathcal{N}(P_c)$
is the regular subdivision induced by the vector $c$.
\end{proof}

We recall the following definition 
from integer programming.  A good reference for these topics is Chapter 22 of Schrijver's book \cite{Sch}.

\begin{definition} \rm
A system of rational inequalities $Dx \leq d$
is called {\it totally dual integral} (TDI) if for each $w \in \Z^m$
such that the linear program $\text{max} \{w \cdot x \, : \, Dx \leq d\}$
has a finite optimal solution, the dual linear program
$\text{min} \{y \cdot d \, : \, yD = d, \,\, y \geq 0 \}$ has an integral
solution.
\end{definition}

The property of being TDI is a property of the given representation
of a polyhedron in terms of inequalities, and not of the polyhedron
itself. In what follows, whenever we say ``the polyhedron $P_c$ is TDI'',
what we mean is that the inequality system 
$\,b_i \cdot x \leq c_i, \,\,\, i=1,\ldots , n$ is TDI.
The following characterization of unimodular configurations
is easily derived from the basic properties of
TDI systems \cite[\S 22]{Sch}.

\begin{proposition}
The vector configuration $B  = \{b_1,\ldots,b_n\} \subseteq \Z^m$
is unimodular if and only if the polyhedron $P_c$
is TDI for every $c \in \Z^n$.
\end{proposition}

We will prove an analogous result for supernormal configurations
by considering only those polyhedra $P_c$ where $c$ ranges
over a certain subset of $\Z^n$.
First we give a name to these special polyhedra.

\begin{definition} \rm The system of inequalities defining $P_c$
is {\it tight} 
if $\,P_{c-e_i}\,\cap \,\Z^m \,$ is  strictly contained in
$\,P_{c}\,\cap \,\Z^m \,$ for every
unit vector $e_i \in \Z^n$.
\end{definition}

Tightness is a property not of the polyhedron $P_c$ but of the
inequality system $\, b_i \cdot x \leq c_i, \,\,\, i=1,\ldots , n $.
However, as with TDI, we shall abuse language by simply saying ``$P_c
$ is tight''.  With this convention, $P_c$ is tight if and only if,
for each $i = 1,\ldots,n$, there exists a lattice point $x \in P_c$
with $b_i \cdot x = c_i$.

\begin{theorem} \label{tight=TDI}
The vector configuration $B  = \{b_1,\ldots,b_n\} \subseteq \Z^m$
is supernormal if and only if every tight polyhedron
$P_c$ is TDI.
\end{theorem}

\begin{proof}
We first prove the if direction using condition (3) in Proposition
\ref{uni-triang}.  Let $\Delta$ be a regular triangulation of $B$
which uses all vectors.  We wish to show that $\Delta$ is a unimodular
triangulation.  By Lemma \ref{regular=fan} there is a simple
polyhedron $P_c$ whose normal fan equals $\Delta$. In particular,
every vector $b_i$ defines a facet of $P_c$.  Since $P_c$ is a
rational polyhedron, there is some $r>0$ such that $P_{rc} = r
P_c$ is integral.  The polyhedron $P_{rc}$ has normal fan
$\Delta$, and is tight, and so is TDI by assumption.  Theorem 22.5 in
\cite{Sch} implies that every set $\sigma$ of $m$ vectors in $B$ that
define a vertex of $P_{rc}$ is a basis of $\Z^m$. These cones
$\sigma$ are the maximal cells of $\Delta$.  Hence $\Delta$ is
unimodular and we conclude that $B$ is supernormal.

For the only-if direction, suppose that $B$ is supernormal and let $c
\in \Z^n$ be such that $P_c$ is tight. Consider any face $F$ of $P_c$,
and let $\sigma$ be the set of all vectors $b_i \in B$ such $\, b_i
\cdot x = c_i\,$ holds for all $x \in F$.  In view of \cite[Theorem
22.5]{Sch}, it suffices to prove that $\sigma$ is a Hilbert basis.
Suppose this is not true.  Supernormality implies
that $\cone(\sigma)$ contains at least one other vector $b_j \in B
\backslash \sigma$.  Because $b_j$ lies in $\cone(\sigma)$, we can
write $b_j=\sum_{b_i \in \sigma} \lambda_i b_i$ where $\lambda_i \geq 0$.
Since $P_c$ is tight
there exists a lattice point $z \in P_c$ with $b_j \cdot z = c_j$.
However since $j \not \in \sigma$, we know that there is some $x \in
F$ for which $b_j \cdot x < c_j$. The first of these two
statements implies $\,c_j=b_j \cdot z = \sum_{b_i \in \sigma}
\lambda_i (b_i \cdot z) \leq \sum_{b_i \in \sigma} \lambda_i c_i$.
The second implies $\,c_j > b_j \cdot x = \sum_{b_i \in
\sigma} \lambda_i (b_i \cdot x) = \sum_{b_i \in \sigma} \lambda_i
c_i$.  But these two statements contradict each other, and so we
conclude that $b_j$ does not exist, and thus $\sigma$ is a Hilbert basis.
It follows that $P_c$ is TDI.
\end{proof}

\section{Chambers and Initial Ideals}
\label{chambers}
The goal of this section is to prove Theorem \ref{gb=chamber}
which states that the chamber complex equals the
Gr\"obner fan if $B$ is supernormal. We start out by
characterizing these two fans by means
of the polyhedra $P_c$. In the next two lemmas,
$B$ is an arbitrary configuration in $\Z^m$.

\begin{lemma} \label{ch-complex}
The chamber complex of $B = \{b_1,\ldots,b_n\} \subseteq \Z^m$  is
the common refinement of the normal fans $\,\mathcal{N}(P_c) \, $ as
$c$ runs over $ \Z^n$.
\end{lemma}

\begin{proof} According to the definition given in the introduction,
two vectors lie in the same cell of the chamber complex
if and only if they lie in exactly the same cones spanned by
linearly independent $m$-subsets of $B$. This holds if and only
if, for every regular subdivision $\Delta$ of $B$, they lie in the
same cell of $\Delta$. Lemma \ref{regular=fan} completes the proof.
\end{proof}

Lemma \ref{ch-complex} coincides with the first statement
in \cite[Proposition 3.3.5]{SST}. The term {\it secondary fan} is 
often used for
the chamber complex. For $c \in \Z^n$ consider the  
lattice polyhedron
$$ Q_c \,\, = \,\, {\rm conv}\, \bigl\{ \, x \in \Z^m \,\,
 : \,\, b_i \cdot x \leq c_i
\,\text{ for } \, i = 1,\ldots,n \,\bigr\}. $$
This is the convex hull of all lattice 
points in the polyhedron $P_c $.

\begin{lemma} \label{GBfan}
The Gr\"obner fan  of the binomial ideal 
$ J_B \subseteq k[x_1,\ldots,x_n] $
 is the common refinement
of the normal fans $\mathcal{N}(Q_c)  $ as
$c$ runs over $ \Z^n$.
\end{lemma}

\begin{proof}
This is the second statement of
\cite[Proposition 3.3.5]{SST}.
\end{proof}

The recipe in the introduction (following equation (\ref{Ubino})) shows
how to derive the initial ideal $in_w(J_B)$ associated with a vector
$w \in \cone(B)$. Note the following subtlety in
our notation: while
$w$ is a vector with $m$ coordinates, it
specifies a term order on monomials in $n$ variables.  

Since $P_c$ is a rational polyhedron there is a positive integer
$r$ such that $rP_c = P_{rc}$ has integer vertices. Hence
$\mathcal{N}(P_c) = \mathcal{N}(P_{rc}) = \mathcal{N}(Q_{rc})$.
This proves the following well-known result:

\begin{corollary} For any configuration
$B = \{b_1,\ldots,b_n\} \subseteq \Z^m$,
the Gr\"obner fan of the ideal $J_B$
refines the chamber complex  of $B$. \qed
\end{corollary}

This says that the cones in the chamber complex of $B$ can split into
smaller cones as one passes to the Gr\"obner fan of $J_B$. It is known
that no splitting happens when $B$ is a unimodular configuration; see
for example \cite[Proposition 8.15 (a)]{GB+CP}.  Theorem
\ref{gb=chamber} says that no splitting happens even when $B$ is only
supernormal. To prove this we need one more lemma:

\begin{lemma}\cite[Corollary 22.1c]{Sch} \label{ifTDI}
$\, \, $ If  $\,P_c \,$ is TDI then $\, P_c = Q_c$.
\qed
\end{lemma}

\smallskip

\noindent {\sl Proof of Theorem \ref{gb=chamber}: } 
Let $B$ be supernormal.  In view of Lemmas \ref{ch-complex} and
\ref{GBfan}, it suffices to prove the following statement: for any $c
\in \Z^n$ there exists $c' \in \Z^n$ such that the normal fan
$\mathcal{N}(Q_c)$ of the integral polyhedron $Q_c$ equals the normal
fan $\mathcal{N}(P_{c'})$ of the rational polyhedron $P_{c'}$.  This
is done by ``pushing in'' all facets of $P_c$ that do not contain
integral points. More precisely, given $c \in \Z^n$, let $x^u$ be the
common divisor of all monomials $x^{c-Bz}$ for $z \in Q_c$.  If
$x^u=1$, then $P_c$ is tight.  Otherwise, $P_{c-u}$ is tight and $Q_c
= Q_{c-u}$. Set $c' = c - u$.  Since $P_{c'}$ is tight, we have that
$P_{c'}$ is TDI by Theorem \ref{tight=TDI}. Using Lemma \ref{ifTDI},
we conclude that $P_{c'} = Q_{c'} = Q_{c}\,$ and hence
$\,\mathcal{N}(Q_c) = \mathcal{N}(P_{c'})$.  \qed

\section{How to subdivide a polygon} \label{polygon}

Let $P$ be a planar convex polygon with integral vertices.
In this section we study  convex vector configurations of the following form:
$$  B \quad = \quad \bigl\{\,
(1,u,v) \in \Z^3 \,\, : \,\, (u,v) \in P \cap \Z^2 \,\bigr\} .  $$
We first show that they are all supernormal.

\begin{proposition} \label{polygonsupernormal}
Every  convex configuration in $\Z^3$ is supernormal.
\end{proposition}

\begin{proof}
Let $B$ be a convex configuration in $\Z^3$ and consider any triangulation
$\Delta$ of $B$ that uses all vectors.  Now a lattice triangle
in the plane which contains no other lattice point has area one half
(by Pick's theorem, for example).  This implies that the triangulation
$\Delta$ is unimodular, and so Proposition \ref{uni-triang} implies
that $B$ is supernormal.
\end{proof}

The {\it chamber complex} of the polygon $P$ is the common refinement
of all lattice triangulations of $P$. Hence the chamber complex of the
vector configuration $B$ is simply the cone over the chamber complex
of $P$.  We draw the chamber complex of $P$ by connecting any pair of
lattice points in $P$ by a straight line segment.

\begin{figure}
\epsfig{file=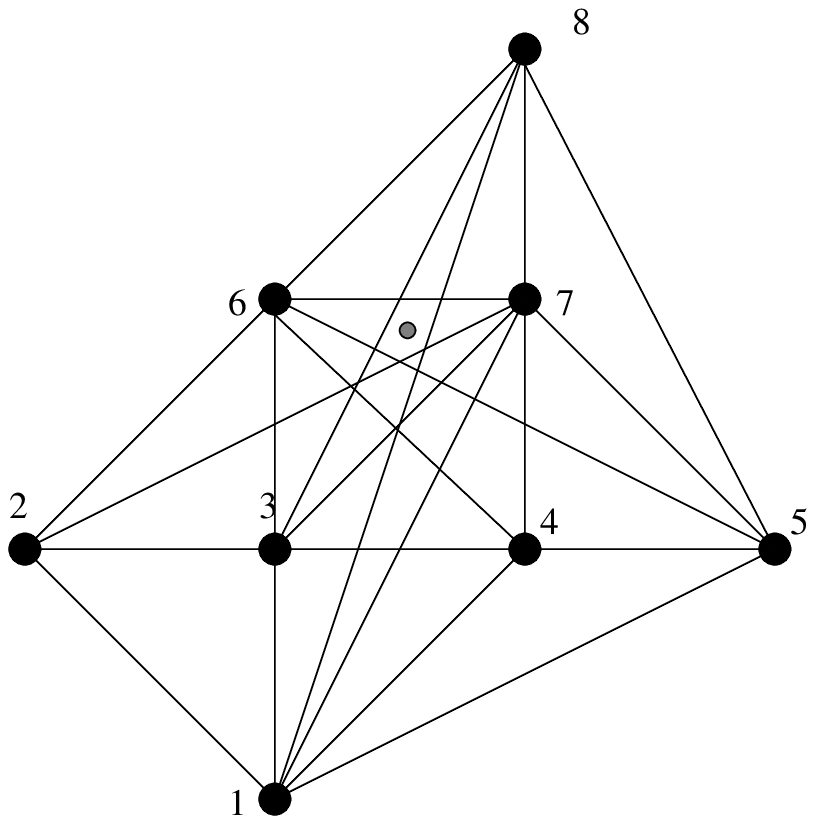, height={8cm}}
\caption{Chamber complex with a pentagonal chamber}
\label{fiveeg}
\end{figure}

For example, if $P$ is the quadrangle with vertices
$\, (1,0), \, (0,1), \, (2,3)$, and $ \, (3,1) $, then $B$
is the set of column vectors of the $3 \times 8$-matrix:
\begin{equation}
\label{threeByeight}
\left(
\begin{array}{llllllll}
1 & 1 & 1 & 1 & 1 & 1 & 1 & 1 \\
1 & 0 & 1 & 2 & 3 & 1 & 2 & 2 \\
0 & 1 & 1 & 1 & 1 & 2 & 2 & 3 \\
\end{array}
\right).
\end{equation}
The chamber complex of $P$ is the subdivision of $P$ into
$26$ triangles, five quadrilaterals, and one pentagon,
which is depicted in Figure \ref{fiveeg}.

We write  $\mu(P)$ for the
maximum number of edges of any region in the chamber
complex of a lattice polygon $P$.
For instance, in  Figure \ref{fiveeg} we have $\mu(P) = 5$.
The main point of this section is the open question of
whether there exists a global upper bound
for the numbers $\mu(P)$.

\begin{problem} \label{PP} {\rm (The Polygon Problem) }
Does there exists a constant $N$ such that
every convex lattice polygon $P$ satisfies $\,\mu(P) \leq N $~?
\end{problem}

We circulated this problem in October 2000, and in the meantime
considerable progress has been made by several people.  However, the
problem remains open for now.\footnote{Two days after we submitted
this paper to a journal, Tracy Hall announced a complete solution to
the Polygon Problem. He constructs a sequence of polygons $P_i$ with
$\,\lim_{i \rightarrow \infty} \mu(P_i) \rightarrow \infty$.}
  Later in this section
we will summarize what is known at the present time (April 2001).

The Polygon Problem is important to us because
it is a special case of a conjecture
in the algebraic theory of integer programming.
Sturmfels and Thomas \cite[Conjecture 6.1]{ST}
asked whether there exists a finite bound $\phi(m)$ on
the number of facets of any cone in the
Gr\"obner fan of an ideal $J_B$ having codimension $m$.
Such a bound would have implications for the
sensitivity analysis of integer programming
in fixed dimension $m$. It is obvious that
$\phi(2) = 2$, and it was conjectured in
\cite[Conjecture 6.2]{ST} that $\phi(3) = 4$.
The latter conjecture  was much too optimistic.
It is now easily seen to be false:
Figure \ref{fiveeg} together with
the following proposition implies $\phi(3) \geq 5$:

\begin{proposition}
Every lattice polygon $P$ satisfies
$\,\phi(3) \geq \mu(P)$.
\end{proposition}

\begin{proof}
The chamber complex of a supernormal
configuration  $B$ is the Gr\"obner fan
of the associated  binomial ideal $J_B$.
Hence $\phi(m)$ is greater or equal to
the maximum number of facets of any cone in the
chamber complex of a supernormal configuration
in $\Z^m$. For $m = 3$ we can take
the chamber complex of a
polygon $P$ to get a lower bound for $\phi(3)$.
\end{proof}

The first counterexamples to \cite[Conjecture 6.2]{ST} were given by
Ho\c{s}ten and Maclagan \cite{HoM} who showed that $\phi(3) \geq
6$. However, the question of  whether $\phi(3)$ is finite remains open.  A
negative answer to the Polygon Problem would show that
 $\phi(m)$ is infinite for  $m \geq 3$.

To  illustrate our algebraic interpretation
of planar chamber complexes, we translate
the marked pentagonal chamber in Figure \ref{fiveeg}
into a specific reduced Gr\"obner basis of binomials.
Our ideal is generated by the three binomials corresponding
to the rows of the matrix in
(\ref{threeByeight}):
$$ J_B
=\langle
x_1x_2x_3x_4x_5x_6x_7x_8-1,x_1x_3x_4^2x_5^3x_6x_7^2x_8^2-1,x_2x_3x_4x_5x_6^2x_7^2x_8^3-1
\rangle.$$
We next fix a term order 
which refines any non-negative real  weight vector $(u_1,\ldots,u_8)$
with the property that $\,w = \sum_{i=1}^8 u_i b_i \,$ lies
in the marked pentagonal chamber of $B = \{b_1,\ldots,b_8\} \in \Z^3$.
For instance, we can take $u = (0, \, 0, \, 0, \, 0, \, 1, \, 4, \, 1, \, 0)$.
The reduced Gr\"obner basis of $J_B$ with respect to this
term order equals:
$$
\{\underline{x_2^4x_3x_6^2-x_4^2x_5^5x_7}, \,\, x_5x_7x_8^2-x_1^2x_2^2x_3,
\,\,
x_2x_6x_8-x_1x_4x_5^2, \,\, \underline{x_7x_8^3-x_1^4x_2^3x_3^2x_4}, $$
$$\underline{x_6x_8^2-x_1^3x_3x_4^2x_5^3},  \,\, x_4x_5^2x_7x_8-x_2,
\,\, x_1x_2^2x_3x_6-x_5, \,\,
\underline{x_1x_4^2x_5^4x_7-x_2^2x_6}, $$
$$\underline{x_1^2x_2x_3x_4x_5-x_8}, \,\,
x_1^2x_3x_4^2x_5^3x_7-1\}
$$
The five ``flippable'' Gr\"obner basis elements
are underlined. They correspond to the five
edges of the pentagonal chamber in 
Figure  \ref{fiveeg}.

\vskip .3cm

We shall now present what is known on the Polygon Problem.
The following result is an outgrowth of the combined efforts
of Miguel Azaola, Jes\'us de Loera, J\"org Rambau, 
Francisco Santos, Marc Pfetsch and G\"unter Ziegler.
In November 2000, the first four of these obtained the lower bound of $12$.
It is attained by the $8 \times 84$ lattice rectangle. In April 2001, the last
two succeeded in improving the previous world record
from $12$ to $15$.  This is the currently best known bound.

\begin{proposition} \label{world-record}
If $P$ is the $9 \times 265$ lattice rectangle
then $\mu(P) = 15$, and  hence $\,
\phi(3) \geq 15$.
\end{proposition}

Pfetsch and Ziegler have made extensive
calculations of the numbers $\mu(P)$ for various lattice rectangles $P$.
Their computational results are posted at the website

\begin{center} 
\url{http://www.math.TU-Berlin.de/~pfetsch/chambers/}
\end{center}

\noindent The data posted at this website  seem to suggest
that the answer to the question in Problem \ref{PP} is 
more likely to be negative.

The example referred to in the proposition above consists
of all lattice points  $(i,j)$ where
$0 \leq i \leq 9$ and $0 \leq j \leq 265$.
Pfetsch and Ziegler identified two chambers which are
$15$-gons in the unit square with 
vertices $(0,132), \,\, (0,133), \,\, (1,132), \,\, (1,133)$.
Note  that one of the edges of this square 
lies on the boundary of the $9 \times 265$
lattice rectangle. It seems that this is not a coincidence:
Ernest Croot has shown that any chamber with many edges must be
located close to the boundary of $P$. The paper with Croot's
precise results is forthcoming.

\section{Virtual Chambers and Virtual Initial Ideals}
\label{virtual}

In Section \ref{chambers} 
we established the bijection 
between chambers of a supernormal configuration $B$ and initial
monomial ideals of $J_B$. In this section we will extend it to a
bijection between {\it virtual} chambers of $B$ and {\it virtual}
initial ideals of $J_B$,  proving Theorem \ref{unimodular2}.
First we define these objects and explain how the bijection works.

Throughout this section we assume that $B = \{b_1, \ldots, b_n\} $ 
generates the lattice $ \Z^m$.  This holds if $B$ is supernormal
by Proposition \ref{uni-triang}.
Under this hypothesis we can find a
configuration $A = \{a_1, \ldots, a_{n}\} \subseteq \Z^{n-m}$ such
that the integer kernel of the $(n-m) \times n$ matrix
$(a_1,\ldots,a_n)$ is spanned by the rows of the matrix
$(b_1,\ldots,b_n)$.  We will also use the 
notation $A$ for the first matrix and $B$ for the second one. 
The relationship between $A$ and $B$ is called
{\it Gale duality} \cite[Chapter 6]{Zi}.  It is well-known
(\cite{BFS}, \cite{BGS}) that the poset of regular subdivisions of $A$
(ordered by refinement) is antiisomorphic to the face poset of the
chamber complex of $B$.

The minimal elements of the poset of regular subdivisions of $A$ are
the regular triangulations of $A$ and they correspond to the
full-dimensional chambers of $B$. This correspondence can be described
explicitly. Let $\Delta = \{\sigma_1, \ldots, \sigma_k\}$ be the
maximal cells of a regular triangulation of $A$ where $\sigma_i =
\{a_{i_1}, \ldots, a_{i_{n-m}}\}$. This defines the chamber
$\bigcap_{t=1}^k \cone({\bar \sigma}_t)$ where ${\bar \sigma}_i = \{
b_j: \,\, j \notin \{i_1, \ldots, i_{n-m}\} \}$.  The bijection
between the regular triangulations of $A$ and the maximal chambers of
$B$ was extended in \cite{DHSS} to {\it all} triangulations of $A$.

\begin{definition} \rm Let $\Delta = \{\sigma_1, \ldots, \sigma_k\}$ be
any (not necessarily regular) triangulation of
the configuration $A$. Then the collection
of  complementary subsets $\{{\bar \sigma}_1, \ldots, {\bar \sigma}_k\}$ of
$B$ is called a {\it virtual chamber} of $B$.
\end{definition}

The configuration in Figure \ref{2x3grid} 
of the Introduction is given by the columns of the
matrix
\begin{equation}
\label{threeBysix}
\left(
\begin{array}{llllll}
1 & 1 & 1 & 1 & 1 & 1 \\
0 & 1 & 2 & 0 & 1 & 2 \\
0 & 0 & 0 & 1 & 1 & 1 \\
\end{array}
\right).
\end{equation}
This configuration $B$ has $18$ virtual chambers.  $16$ of these are
chambers and hence visible in Figure \ref{2x3grid}. The two additional
virtual chambers are
$$\{
(1, 3, 4),  
(1, 3, 5), 
(1, 4, 6),
(1, 5, 6),  
(2, 3, 4),  
(2, 3, 5),  
(2, 4, 6), 
(2, 5, 6)\},$$
$$\{
(1, 2, 5),  
(1, 2, 6),
(1, 3, 5),  
(1, 3, 6),  
(2, 4, 5),  
(2, 4, 6),  
(3, 4, 5),  
(3, 4, 6)
\}.$$
We invite the reader to ``locate'' these virtual chambers
in Figure \ref{2x3grid}, and to draw the 
two non-regular triangulations of
$A$. 

We define an $(n-m)$-dimensional grading of the polynomial ring 
$S = k[x_1, \ldots, x_n]$ by setting the degree of  $x_i$ 
to be $a_i$ for $i=1, \ldots, n$. Thus $S$ is graded by the
monoid $\N A$ which is spanned by the Gale dual  configuration $A$.
 The ideal $J_B$ is homogeneous in
this grading since
$$x^u - x^v \in J_B \quad
\text{  if and only if  } \quad \sum_{i=1}^n u_ia_i
= \sum_{i=1}^n v_ia_i.$$
The Hilbert function of the quotient ring $S / J_B$ is given by
\begin{equation}
\label{HilF}
\dim_k((S / J_B)_b) \quad =
\quad  \bigg\lbrace \aligned 1 & \text{ if } 
b \in \N A, \\             0 & \text{ otherwise.}
\endaligned
\end{equation}
A homogeneous ideal in $S$ with the same Hilbert function as $J_B$ was
called an $A$-graded ideal in \cite[\S 10]{GB+CP}.
Monomial $A$-graded ideals include, but are not limited to, initial
ideals of $in_w(J_B)$.

\begin{definition} \label{DefofVirtual}
 \rm A monomial ideal $M$ in $S$ is a {\it virtual initial ideal} of
 $J_B$ if the Hilbert function of $S / M$ is equal to the Hilbert
 function (\ref{HilF}).  This means that for every degree 
$b \in \N A$ there is  exactly one monomial $x^u$ 
of degree $b$ with the property that $x^u  \not \in M$.
\end{definition}

To illustrate this definition and 
 Theorem \ref{unimodular2} we compute a
virtual initial ideal of $J_B$ for 
(\ref{threeBysix}). First consider $w = (2,2,1)$. Then
$$ in_w(J_B) = 
\langle x_1 x_2 x_3, x_4 x_5 x_6, x_3 x_5 x_6^2,
x_1^2 x_2 - x_5 x_6^2,
x_4^2 x_5 - x_2 x_3^2,
x_3 x_6 - x_1 x_4
 \rangle $$
This $A$-graded ideal corresponds to the centroid in Figure \ref{2x3grid}.
By replacing each of the three binomials by one of its terms, we
get eight virtual initial ideals of $J_B$, one for each
virtual chamber adjacent to the centroid in Figure \ref{2x3grid}.  
For instance, taking the first term in each of the
three binomial generators of $in_e(J_B)$ gives
the virtual initial ideal
\begin{eqnarray*} &
 \langle x_1^2, x_3, x_4\rangle \, \cap \,
 \langle x_1^2, x_3, x_5\rangle \, \cap \,
 \langle x_1, x_4^2, x_6\rangle \, \cap  \,
 \langle x_1, x_5, x_6\rangle \, \cap
 \\ &
\!\! \langle x_2, x_3, x_4\rangle  \cap 
 \langle x_2, x_3, x_5\rangle  \cap 
 \langle x_2, x_4^2, x_6\rangle \cap 
 \langle x_2, x_5, x_6\rangle  \cap
 \langle x_1^2, x_3, x_4^2,x_6\rangle.
\end{eqnarray*}
We pass to the radical of this ideal by erasing all exponents, and
deleting the embedded component at the end.  The eight remaining index
sets are precisely the cells in the first virtual chamber listed for
example (\ref{threeBysix}).  This process of using primary
decomposition to read off the virtual chamber from a given virtual
initial ideal works in general:

\begin{remark}
The map referred to in Theorem \ref{unimodular2}
is given by:
\begin{equation} \label{Mappp}
M \,\,\, \mapsto \,\,\,
\bigl\{ \, \bar \sigma \,\, \, : \,\,\,
\langle x_i \, : \, i \in  \bar \sigma \rangle
\text{ is a minimal prime of } \,  M \, \bigr\} .
\end{equation}
\end{remark}

This remark  follows  essentially from \cite[Theorem 10.10]{GB+CP}.
We shall give an alternative description
of the map (\ref{Mappp}) after Lemma \ref{Section} below.
That description will be self-contained, 
with no reference to \cite{GB+CP} needed,
and better suited for the purpose of proving
Theorem \ref{unimodular2}.

For arbitrary configurations $B$, the map (\ref{Mappp}) is neither
injective nor surjective. Two virtual initial ideals can give rise to
the same virtual chamber and there might be virtual chambers which do
not correspond to virtual initial ideals \cite[Theorem 10.13]{GB+CP}.
What we are claiming in Theorem \ref{unimodular2} is that for
supernormal configurations $B$ the map (\ref{Mappp}) is both
injective and surjective.  In the special case when $B$ is unimodular
this was proved in \cite[Lemma 10.14]{GB+CP}.

We next present a characterization of
virtual monomial ideals in terms of the
integral polyhedra $Q_c$ introduced in Section 4.

\begin{lemma} \label{virtual-fiber}
A monomial ideal $M$ is a virtual initial ideal
of $J_B$ if and only if, for every  $c \in \Z^n$,
the polyhedron $Q_c$ is either empty or $Q_c$ contains
a unique lattice point $z$ such that 
$\, \prod_{i=1}^n x_i^{c_i - b_i \cdot z}  \,$ is not in $M$.
\end{lemma}

\begin{proof} The map $\, z \mapsto  \prod_{i=1}^n x_i^{c_i - b_i \cdot z} \,$
is a bijection between the  set of lattice points in $Q_c$ 
and the set of monomials in $S$ having  degree $\, \sum_{i=1}^n c_i a_i $.
Hence the condition in the lemma states that
every non-zero graded component of $S$ contains
exactly one monomial which is not in $M$.
\end{proof}

In \cite[Proposition 10.8]{GB+CP} it was shown that the lattice point
$z$ chosen as in Lemma \ref{virtual-fiber} 
need not be a vertex of the polyhedron $Q_c$.
This is not the case for initial ideals of $J_B$, 
and the following
important lemma states that it is also not the case if $Q_c=P_c$.

\begin{lemma} \label{vvvertex}
If $P_c$ is non-empty and equal to $\,Q_c = {\rm conv}( P_c \,\cap
\,\Z^m)\,$ then the lattice point $z$ selected in Lemma
\ref{virtual-fiber} is a vertex of $Q_c$.
\end{lemma}

\begin{proof}
Let $z_1, \ldots, z_r$ be the vertices of $\, P_c = Q_c \,$ and let
$x^{u_1}, \ldots, x^{u_r}$ be the corresponding monomials in $S$ of
degree $\,b = \sum_{i=1}^n c_i a_i $.  

We first show that every monomial in $S_{r b}$ lies in the monomial
ideal $\, \langle x^{u_1}, \ldots, x^{u_r} \rangle \, \subseteq \, S$.
In polyhedral terms, if $z$ is any lattice point in $P_{r c} =
Q_{r c}$, then $z$ can be written as $\, z
\, = \, \sum_{i=1}^r \gamma_i z_i \, + \, w \,$ where $w \in P_0 $ and
the $\gamma_i$ are non-negative reals summing to $r$.  This means
that for the corresponding monomial $x^u$ we have $u = r c-Bz= r
c- \sum_{i=1}^r \gamma_i Bz_i - Bw = \sum_{i=1}^r \gamma_i u_i - Bw$, 
where $Bw \in (\mathbb Z_{\leq 0})^n$, since $w \in P_0$.  
There exists an index $j \in \{1,\ldots,r\}$ such that
$\,\gamma_j \geq 1\,$ and this implies that  $u \geq u_j$, and thus
that $x^{u_j} $ divides $ x^u$.  This shows that $S_{rb}$ lies in 
$\, \langle x^{u_1}, \ldots, x^{u_r} \rangle $.

Since our virtual initial ideal $M$ must have a standard monomial of
degree $r b$, it cannot contain the ideal $\, \langle x^{u_1},
\ldots, x^{u_r} \rangle $, and we conclude that one of the monomials
$\,x^{u_j} \,$ is not in $M$, as desired.
\end{proof}

\smallskip

We next present an alternative characterization of triangulations of
$A$, and hence of virtual chambers of $B$.  A subset $U$ of the closed
orthant $\R^n_+$ is an {\it order ideal} if $v \in U$ and $ u \leq v $
coordinatewise implies $u \in U$.  Let $\pi$ be the linear map
$\,(\lambda_1 , \ldots, \lambda_n) \mapsto \sum_{i=1}^n \lambda_i a_i
\,$ from $\R^n_+$ onto ${\rm cone}(A)$.  A section of $\pi$ is a
map $\, s : {\rm cone}( A) \rightarrow \R^n_+\,$ such that the
composition $\,\pi \circ s \,$ is the identity on ${\rm
cone}(A)$. Note that every triangulation $\,\Delta \,$ of $A$ defines
a section $s_\Delta$ as follows: $s_\Delta(b)$ is the unique
vector $u \in \R_+^n$ with $Au = b$ and whose support is a cell of
$\Delta$. The image $\,{\rm im}(s_\Delta)\,$ 
of such a section $\,s_\Delta \,$ is an order
ideal in $\R^n_+$.

\begin{lemma} \label{Section}
The map $\Delta \mapsto s_\Delta$ is a bijection between
triangulations of $A $ and sections $s$ of $\pi$ for which
${\rm im}(s)$ is an order ideal in $\R_+^n$.
\end{lemma}

\begin{proof}
It is clear that the section $s_{\Delta}$ associated to a
triangulation $\Delta$ of $A$ satisfies the desired conditions, so we
need only show that every section $s$ satisfying the hypothesis 
comes from a triangulation.

Fix such an $s$.
We first observe that $s(rb)=rs(b)$ for
$b \in {\rm cone}(A)$ and $r \in \R_+$.  If
$r<1$ then $c=rs(b) \in {\rm im}(s)$, and so
$\pi(c)=r\pi(s(b))=rb$ satisfies $s(rb)=c=rs(b)$.  The case that
$r>1$ follows from this.

We claim that the set of all possible supports of vectors in 
${\rm im}(s)$ is a
triangulation of $A$.  We first show that the subsets of $A$
indexed by these supports are linearly independent.
  Suppose not, so for some $b \in \mathbb R^n_+$ there is a vector 
$u = (u_1, \ldots, u_n)$ such that $Au = b$ 
where ${\rm supp}(u)$ is a proper subset of ${\rm supp}(s(b))$.  
There is some $r>0$ for which $r u < s(b)$, and so
$ru \in {\rm im}(s)$.  Now $\pi(ru)=r \pi(u)=rb$, so $s(rb)=ru.$ 
This implies that $s(b)=u$, a contradiction
since ${\rm supp}(s(b))$  properly contains ${\rm supp}(u)$.

This shows that the cones
$\cone(a_i: i \in {\rm supp}(s(b)))$ as $b$ ranges over $\cone(A)$ 
are simplicial and that they  cover $\cone(A)$.  We also note that this
argument actually shows that for any $b'$ in the 
relative interior of
$ \cone( a_i : i \in {\rm supp}(s(b)))$ we have ${\rm supp}(s(b'))=
{\rm supp}(s(b))$.  Hence  the relative interiors of two distinct cones do not
intersect.  The order ideal hypothesis guarantees that these cones form
a simplicial fan. 
\end{proof}

This bijection means we can express the map in Theorem 
\ref{unimodular2} as taking a virtual
initial ideal $M$ to a section $s$ such that $\, {\rm im}(s)\,$ is
an order ideal in $\R^n_+$.  Fix $M$.  
For $P_c = Q_c$ \, we set $s(Ac) = c - Bz$ where $z$ is given by
Lemma \ref{vvvertex}. Since $s(rb) = rs(b)$ we can extend this to 
all rational $P_c$, and hence to all $\, b \in \cone(A) \,$ by continuity. 

Now we are ready to prove  Theorem \ref{unimodular2}.  Recall that a
polyhedron $P_c$ is tight if and only if the greatest common divisor
of all monomials of the form $x^{c-Bz}$ for $z \in P_c$ is
one.  If $P_c$ is not tight, let $x^w$ be the greatest common
divisor of all monomials of the corresponding degree. Then if a 
monomial ideal $I$ is generated in tight degrees, $x^u \not \in I$
implies $x^{u-w} \not \in I$ where $u=c-Bz$ for some $z \in P_c$.
We first present the part of the proof that holds for a general
configuration.

\begin{lemma} \label{nontight}
Let $x^u$ divide $x^v$, and let $x^w$ and $x^{w'}$ be the greatest common
divisors of all monomials of the same degree as $x^u$ and $x^v$
respectively.  Then $x^{u-w}$ divides $x^{v-w'}$.
\end{lemma}

\begin{proof}

Suppose this is not the case, so there is some $i$ with $(u-w)_i >
(v-w')_i$.  Since the greatest common divisor of all monomials of the
same degree as $x^{u-w}$ is $1$, we know that $P_{u-w}$ is tight, and so
there is some lattice point $z \in P_{u-w}$ such that $b_i \cdot z =
(u-w)_i$.  Because $u-w < v$, we also have $z \in P_{v}$.  This means
$x^{v-Bz}$ is a monomial of the same degree as $x^{v}$, and is thus
divisible by $x^{w'}$, so $v-w'-Bz \geq 0$.  But this implies that
$b_i \cdot z =(u-w)_i \leq (v-w')_i$, a contradiction.
\end{proof}

\smallskip

\noindent {\sl Proof of Theorem \ref{unimodular2}}: For each virtual
chamber of $B$ we will construct a virtual initial ideal which maps to
it.  The construction will make it clear that this map is injective.
Let $s$ be the section of $\pi$ corresponding to our virtual
chamber, as described in Lemma \ref{Section}.  It is straightforward
to check that $s(Ac)$ is a vertex of the polyhedron $P_c$
for  every $c \in \R_+^n$. 

We define $M$ to
be the ideal generated by all monomials $x^c$ such that $P_c$ is tight
and $c$ is not in the image of $s$.  We claim that $M$ is a virtual
initial ideal.
By construction, $M$ has at most one standard monomial in every tight
degree, and thus in every degree. Tight polyhedra are 
integral
by Theorem \ref{tight=TDI} and Lemma \ref{ifTDI}. If
$P_c$ is tight then $s(Ac)$ is a vertex of $P_c = Q_c$ and hence
$s(Ac) \in \N^n$.
We claim $x^{s(Ac)} \not \in M$
for all $c$ such that $P_c$ is tight.  If not, there is some generator
$x^v$ of $M$ with $P_v$ tight dividing $x^{s(Ac)}$.  
But since ${\rm im}(s)$ is an order
ideal, we must have $v \in {\rm im}(s)$, contradicting $x^v \in M$.  
Therefore $x^{s(Ac)} \not \in M$.

If $P_c$ is not tight, let $x^w$ be the greatest common divisor
of all monomials of degree $Ac$.
Then we claim that $x^{u+w} \not \in M$, where
$x^u \not \in M$ satisfies $u=c-w-Bz$ for $z \in P_{c-w}$.  Otherwise
 there would be some generator $x^v$ of $M$
with $x^v$ dividing $x^{u+w}$.  But since $P_v$ would then be a tight
degree, Lemma \ref{nontight} would imply that $x^v$ must divide $x^u$,
a contradiction.  This concludes the proof that $M$ is a virtual initial ideal.

The virtual initial ideal $M$ just constructed is clearly mapped back
to $s$ under the map  (described after Lemma \ref{Section})
from virtual initial ideals to triangulations.
Hence this map is a bijection as desired.
\qed

\smallskip 

\noindent
{\bf Acknowledgement.} We thank Miguel Azaola, Ernest Croot, Jes\'us
de Loera, Tracy Hall, Marc Pfetsch, 
J\"org Rambau, Francisco Santos, and G\"unter
Ziegler for permitting  us to discuss their work in Section~\ref{polygon}.

\end{document}